\documentclass[12pt]{article}

\usepackage{graphicx,amssymb,amsmath,amsthm,amsfonts}
\usepackage[top=3cm, bottom=3cm, right=3cm, left=3cm]{geometry}

\begin{document}

\pagenumbering{arabic}
\title{\huge \bf Lie point symmetries and conservation laws for a Gardner type system}
\author{\rm \large Valter Aparecido Silva Junior$$ \\ \\ \\
\it Instituto Federal de Educa\c{c}\~ao, Ci\^encia e \\
\it Tecnologia de S\~ao Paulo - IFSP \\
\it Avenida Marginal, 585, Fazenda Nossa Senhora Aparecida do Jaguari\\ \it   $13871$-$298$ - S\~ao Jo\~ao da Boa Vista/SP, Brasil \\
\rm E-mail: valtersjunior@ifsp.edu.br \\ \\ \\
\it Instituto de F\'\i sica ``Gleb Wataghin" - IFGW \\
\it Universidade Estadual de Campinas - UNICAMP \\
\it Rua S\'ergio Buarque de Holanda, 777, Cidade Universit\'aria\\ \it   $13083$-$859$ - Campinas/SP, Brasil \\
\rm E-mail: valtersjunior@gmail.com}
\date{\ }
\maketitle

\newpage 

\begin{abstract}
We determine the Lie point symmetries of a Gardner type system and establish its nonlinear self-adjointness. We then construct conservation laws via Ibragimov's Theorem.
\end{abstract}

\vskip 1cm

\begin{center}
{2010 AMS Mathematics Classification numbers:\vspace{0.2cm}\\
76M60, 58J70, 70G65

\

}
\end{center}

\newpage

\section{Introduction}

The well-known KdV, mKdV and Gardner equations \cite{fll,lzxs,m}, given respectively by
\begin{equation}
u_t + auu_x + cu_{xxx} = 0, \quad u_t + bu^2u_x + cu_{xxx} = 0\nonumber
\end{equation}
and
\begin{equation}
u_t + (au + bu^2)u_x + cu_{xxx} = 0,\nonumber
\end{equation}
are special cases of
\begin{equation}
\label{ggardner}
u_t + (au^p + bu^{2p})u_x + cu_{xxx} = 0, \quad a, b, c \in \mathbb{R}, \quad p \in \mathbb{Z},
\end{equation}
a recurrent nonlinear evolution equation that, in the framework of hydrodynamics, models internal gravity waves in a density-stratified ocean \cite{bik}. Of great physical importance for also describing a variety of wave phenomena in plasma, solid state and quantum physics, its integrability has been investigated by several authors under different approaches. Explicit solutions (including solitons) and a detailed discussion of their properties (like stability) can be found in \cite{hmhs,mrz,zlbb} and references therein.

The situations covered by the above equations tend to be idealized. This occurs, among some reasons, because such models are derived without considering the possible effects of one wave on another. However, since the pioneer works of Hirota and Satsuma \cite{hs,hs1}, where they present a system to study the interactions of two long waves with distinct dispersion relations, an increasing number of coupled KdV equations have been introduced in the literature with a substantially larger and more realistic range of applications. See for example \cite{s,tjb,ws}. In particular, starting from \cite{aam,ti}
\begin{equation}
\left\{
\begin{aligned}
u_t + a[(u^2 + v^2)u]_x + u_{xxx} = 0&\\
v_t + a[(u^2 + v^2)v]_x + v_{xxx} = 0&
\end{aligned}
\right.
\nonumber
\end{equation}
and
\begin{equation}
\left\{
\begin{aligned}
u_t + (u^pv^{p + 1})_x + u_{xxx} = 0&\\
v_t + (u^{p + 1}v^p)_x + v_{xxx} = 0&
\end{aligned}
\right.,
\nonumber
\end{equation}
we propose in this paper the new class of system
\begin{equation}
\label{pKdV}
\left\{
\begin{aligned}
F_1 \equiv u_t + a[(u^p + v^p)u]_x + b[(uv)^pv]_x + cu_{xxx} = 0&\\
F_2 \equiv v_t + a[(u^p + v^p)v]_x + b[(uv)^pu]_x + cv_{xxx} = 0&
\end{aligned}
\right.,
\end{equation}
a two-component generalization\footnote{If $u = v$, the system (\ref{pKdV}) is reduced to equation (\ref{ggardner}).} of the equation (\ref{ggardner}), with the objective of determining the Lie point symmetries (Section 2). Established its nonlinear self-adjointness (Section 3), we then construct conservation laws via the recent Ibragimov's Theorem  (Section 4), an extension of the celebrated Noether's Theorem to problems with no variational structure. The main references used are \cite{bk,gan,i1,i3,i2,i4,i5,ol}. 

In what follows, we assume that $\{ac, bc\} \neq \{0\}$ and $p \geq 1$. All constants $c_i$'s are arbitrary and all functions are smooth.

\section{Lie Point Symmetries Classification}

To begin with, we denote
\begin{equation}
f(u, v) \equiv a(u^p + v^p)u + b(uv)^pv, \quad g(u, v) \equiv a(u^p + v^p)v + b(uv)^pu.
\nonumber
\end{equation}
Applying the standard algorithm presented in \cite{bk} and \cite{ol}, a differential operator
\begin{equation}
X = \mathcal{T}(t, x, u, v)\frac{\partial}{\partial t} + \mathcal{X}(t, x, u, v)\frac{\partial}{\partial x} + \mathcal{U}(t, x, u, v)\frac{\partial}{\partial u} + \mathcal{V}(t, x, u, v)\frac{\partial}{\partial v}\nonumber
\end{equation}
generates the Lie point symmetries of the system (\ref{pKdV}) if the conditions of invariance
\begin{equation}
\label{cond1}
\begin{aligned}
(&D_tW^u + \mathcal{T}u_{tt} + \mathcal{X}u_{tx}) + f_u(D_xW^u + \mathcal{T}u_{tx} + \mathcal{X}u_{xx}) + f_v(D_xW^v + \mathcal{T}v_{tx} + \mathcal{X}v_{xx}) \ +\\
&+ f_{uu}\mathcal{U}u_x + f_{uv}(\mathcal{V}u_x + \mathcal{U}v_x) + f_{vv}\mathcal{V}v_x + c(D_x^3W^u + \mathcal{T}u_{txxx} + \mathcal{X}u_{xxxx}) = MF_1 + NF_2
\end{aligned}
\end{equation}
and
\begin{equation}
\label{cond2}
\begin{aligned}
(&D_tW^v + \mathcal{T}v_{tt} + \mathcal{X}v_{tx}) + g_u(D_xW^u + \mathcal{T}u_{tx} + \mathcal{X}u_{xx}) + g_v(D_xW^v + \mathcal{T}v_{tx} + \mathcal{X}v_{xx}) \ +\\
&+ g_{uu}\mathcal{U}u_x + g_{uv}(\mathcal{V}u_x + \mathcal{U}v_x) + g_{vv}\mathcal{V}v_x + c(D_x^3W^v + \mathcal{T}v_{txxx} + \mathcal{X}v_{xxxx}) = PF_1 + QF_2,
\end{aligned}
\end{equation}
where
\begin{equation}
W^u = \mathcal{U} - \mathcal{T}u_t - \mathcal{X}u_x, \quad W^v = \mathcal{V} - \mathcal{T}v_t - \mathcal{X}v_x,\nonumber
\end{equation}
are satisfied for some set of coefficients $M$, $N$, $P$ and $Q$ to be determined. From (\ref{cond1}) and (\ref{cond2}), the determining equations can be written as
\begin{equation}
\label{determining-equations}
\begin{gathered}
M = \mathcal{U}_u - \mathcal{T}_t, \ N = \mathcal{U}_v, \ P = \mathcal{V}_u, \ Q = \mathcal{V}_v - \mathcal{T}_t,\\
\mathcal{T}_x = \mathcal{T}_u = \mathcal{T}_v = \mathcal{T}_t - 3\mathcal{X}_x = 0,\\
\mathcal{X}_t = \mathcal{X}_u = \mathcal{X}_v = 0,\\
\mathcal{U}_t = \mathcal{U}_x = \mathcal{U}_{uv} = \mathcal{U} - (\mathcal{U}_uu + \mathcal{U}_vv) = 0,\\
\mathcal{V}_t = \mathcal{V}_x = \mathcal{V}_{uv} = \mathcal{V} - (\mathcal{V}_uu + \mathcal{V}_vv) = 0,\\
a[2\mathcal{X}_x + p\mathcal{U}_u + (2 - p)\mathcal{V}_u] = 0,\\
b(\mathcal{X}_x + p\mathcal{U}_u)= b\mathcal{U}_v = b\mathcal{V}_u = 0,\\
(1 - p)(\mathcal{U}_v + \mathcal{V}_u) = (1 - p)(2 - p)\mathcal{U}_v = 0,\\
\mathcal{U}_u - \mathcal{V}_v = (2 - p)(\mathcal{U}_v - \mathcal{V}_u).
\end{gathered}
\nonumber
\end{equation}
It's easy to see therefore that
\begin{equation}
\begin{gathered}
\mathcal{T} = 3c_1t + c_2, \ \mathcal{X} = c_1x + c_3,\\
\mathcal{U} = c_4v + [c_6 + (2 - p)c_4]u, \ \mathcal{V} = [c_6 + (2 - p)c_5]v + c_5u
\end{gathered}
\nonumber
\end{equation}
with
\begin{equation}
\begin{gathered}
a(2c_1 + c_4 + c_5 + pc_6) = 0,\\
bc_4 = bc_5 = b(c_1 + pc_6) = 0,\\
(1 - p)(2 - p)c_4 = (1 - p)(c_4 + c_5) = 0.
\end{gathered}
\nonumber
\end{equation}

\

\noindent{{\bf Proposition 1.} \textit{The Lie point symmetries of the system $(\ref{pKdV})$ are summarized in Table $1$, where}
\begin{equation}
\begin{gathered}
X_1 = 3t\frac{\partial}{\partial t} + x\frac{\partial}{\partial x} - \frac{1}{p}\!\left(u\frac{\partial}{\partial u} + v\frac{\partial}{\partial v}\right)\!\!,\\
Y_1 = 3t\frac{\partial}{\partial t} + x\frac{\partial}{\partial x} - \frac{2}{p}\!\left(u\frac{\partial}{\partial u} + v\frac{\partial}{\partial v}\right)\!\!, \ X_2 = \frac{\partial}{\partial t}, \ X_3 = \frac{\partial}{\partial x},\\
X_4 = v\!\left(\frac{\partial}{\partial u} - \frac{\partial}{\partial v}\right)\!\!, \ Y_4 = v\frac{\partial}{\partial u} - u\frac{\partial}{\partial v}, \ X_5 = u\!\left(\frac{\partial}{\partial u} - \frac{\partial}{\partial v}\right)\!\!.
\end{gathered}
\nonumber
\end{equation}

\begin{table}[h!]
\centering
\begin{tabular}{ccc}
\hline
\multicolumn{1}{|c|}{$a = 0$}                                                      & \multicolumn{1}{c|}{$b = 0$}                                                                                                   & \multicolumn{1}{c|}{$ab \neq 0$}  \\ \hline
\multicolumn{1}{|c|}{\begin{tabular}[c]{@{}c@{}}$X_1$, $X_2$\\ $X_3$\end{tabular}} & \multicolumn{1}{c|}{\begin{tabular}[c]{@{}c@{}}$Y_1$, $X_2$, $X_3$, $X_4$ ($p = 1$)\\ $Y_4$ ($p = 2$), $X_5$ ($p = 1$)\end{tabular}} & \multicolumn{1}{c|}{$X_2$, $X_3$} \\ \hline
\multicolumn{1}{l}{}                                                               & \multicolumn{1}{l}{}                                                                                                           & \multicolumn{1}{l}{}             
\end{tabular}
\\
\tablename{ 1.}
\end{table}

\section{Self-Adjointness Classification}

For the system (\ref{pKdV}), the formal Lagrangian is
\begin{equation}
\mathcal{L} = \bar{u}F_1 + \bar{v}F_2.\nonumber
\end{equation}
Here $\bar{u}$ and $\bar{v}$ are the new dependent variables (so-called nonlocal variables). Calculated the adjoint equations
\begin{equation}
\label{adj1}
F_1^* \equiv -\frac{\delta\mathcal{L}}{\delta u} = \bar{u}_t + f_u\bar{u}_x + g_u\bar{v}_x + c\bar{u}_{xxx} = 0\nonumber
\end{equation}
and
\begin{equation}
\label{adj2}
F_2^* \equiv -\frac{\delta\mathcal{L}}{\delta v} = \bar{v}_t + f_v\bar{u}_x + g_v\bar{v}_x + c\bar{v}_{xxx} = 0,\nonumber
\end{equation}
where $\delta/\delta u$ and $\delta/\delta v$ are Euler-Lagrange operators, we impose that
\begin{equation}
\label{adjointness-condition}
F_1^*|_{(\bar{u}, \bar{v}) = (\varphi, \psi)} = MF_1 + NF_2, \quad F_2^*|_{(\bar{u}, \bar{v}) = (\varphi, \psi)} = PF_1 + QF_2.
\end{equation}
Again $M$, $N$, $P$ and $Q$ is a set of coefficients to be determined and
\begin{equation}
\label{sub}
\varphi = \varphi(t, x, u, v), \quad \psi = \psi(t, x, u, v)
\end{equation}
two functions that not vanish simultaneously. 

As
\begin{equation}
F_1^*|_{(\bar{u}, \bar{v}) = (\varphi, \psi)} = D_t\varphi + f_uD_x\varphi + g_uD_x\psi + cD_x^3\varphi
\nonumber
\end{equation} 
and
\begin{equation}
F_2^*|_{(\bar{u}, \bar{v}) = (\varphi, \psi)} = D_t\psi + f_vD_x\varphi + g_vD_x\psi + cD_x^3\psi,
\nonumber
\end{equation} 
from (\ref{adjointness-condition}) it's possible to conclude that $M = \varphi_u$, $N = \varphi_v$, $P = \psi_u$, $Q = \psi_v$ and
\begin{equation}
\begin{gathered}
\varphi_t - \psi_t = \varphi_x - \psi_x = \varphi_u - \psi_v = \varphi_v - \psi_u = 0,\\
\varphi_t + 2a(u + v)\varphi_x = b\varphi_x = (p - 1)\varphi_x = 0,\\
a(p - 2)(\varphi_u - \varphi_v) = b\varphi_v = (p - 1)\varphi_v = 0,\\
\varphi_{xx} = \varphi_{xu} = \varphi_{xv} = \varphi_{uu} = \varphi_{uv} = \varphi_{vv} = 0.
\end{gathered}
\nonumber
\end{equation}
Hence
\begin{equation}
\left\{
\begin{aligned}
\varphi = (2ac_1t + c_2)v + (2ac_1t + c_3)u - c_1x + c_4\\
\psi = (2ac_1t + c_3)v + (2ac_1t + c_2)u - c_1x + c_5
\end{aligned}
\right.,
\nonumber
\end{equation}
with
\begin{equation}
\label{cod2}
\begin{gathered}
(p - 1)c_1 = bc_1 = (p - 1)c_2 = bc_2 = 0,\\
a(p - 2)(c_2 - c_3) = 0.
\end{gathered}
\nonumber
\end{equation}

\

\noindent{{\bf Proposition 2.}} \textit{The system $(\ref{pKdV})$ is nonlinearly self-adjoint. The substitutions $(\ref{sub})$ are as follows.}

\

\textbf{i)} If $a = 0$ or $p = 2$,
\begin{equation}
\varphi = c_3u + c_4, \quad \psi = c_3v + c_5;
\nonumber
\end{equation}

\

\textbf{ii)} if $b = 0$ and $p = 1$,
\begin{equation}
\varphi = (2ac_1t + c_2)(v + u) - c_1x + c_4, \quad \psi = (2ac_1t + c_2)(v + u) - c_1x + c_5;\nonumber
\end{equation}

\

\textbf{iii)} in all other cases,
\begin{equation}
\varphi = c_4, \quad \psi = c_5.
\nonumber
\end{equation}

\

\noindent{\textbf{Remark 1.} Actually, the system (\ref{pKdV}) is quasi self-adjoint. In particular, it's strictly self-adjoint iff $a = 0$ or $p = 2$.}

\

\noindent{\textbf{Remark 2.} Indeed, any system of the form}
\begin{equation}
\left\{
\begin{aligned}
u_t + r_x(u, v) + cu_{xxx} = 0\\
v_t + s_x(u, v) + cv_{xxx} = 0
\end{aligned}
\right.
\nonumber
\end{equation}
is strictly self-adjoint whenever $r_v = s_u$.

\section{Conservation Laws}

According to the theorem demonstrated by Ibragimov \cite{i3}, the components of the conserved vector $C = (C^t, C^x)$ associated to $X$, a Lie point symmetry admitted by the system (\ref{pKdV}), are given by
\begin{equation}
\label{comp-t}
C^t = \varphi W^u + \psi W^v
\end{equation}
and
\begin{equation}
\label{comp-x}
\begin{aligned}
C^x = \varphi&(f_uW^u + f_vW^v) + c(\varphi D_x^2 - D_x\varphi D_x + D_x^2\varphi)W^u +\\
&+ \psi(g_uW^u + g_vW^v) + c(\psi D_x^2 - D_x\psi D_x + D_x^2\psi)W^v.
\end{aligned}
\end{equation}

In view of Proposition 2, the next result takes into account all the generators of Table 1. In most cases, however, the expressions (\ref{comp-t}) and (\ref{comp-x}) lead us to trivial conservation laws or to the vectors
\begin{equation}
C^t = u, \quad C^x = f(u, v) + cu_{xx}\nonumber
\end{equation}
and
\begin{equation}
C^t = v, \quad C^x = g(u, v) + cv_{xx}\nonumber
\end{equation}
that can be readily obtained from the system (\ref{pKdV}) by direct integration. Below we only bring the unobvious cases.

\

\noindent{{\bf Proposition 3.} \textbf{i)} {\it Let $a = 0$. For}
\begin{equation}
\varphi = u, \quad \psi = v,\nonumber
\end{equation}
$X_1$ provides}
\begin{equation}
\begin{gathered}
C^t = (p + 1)(u^2 + v^2), \\ C^x = 2b(2p + 1)(uv)^{p + 1} + c(p + 1)[(u^2 + v^2)_{xx} - 3(u_x^2 + v_x^2)].\nonumber
\end{gathered}
\end{equation}

\

\noindent{\textbf{ii)} {\it Let $b = 0$.}}

\

\textbf{ii.a)} {\it Taking}
\begin{equation}
\varphi = \psi = u + v\nonumber
\end{equation}

for $Y_1$ and 
\begin{equation}
\varphi = \psi = 2at(u + v) - x\nonumber
\end{equation}

for $X_2$, we obtain
\begin{equation}
\begin{gathered}
C^t = 3(u + v)^2, \\ C^x = 4a(u + v)^3 + 3c\{[(u + v)^2]_{xx} - 3(u_x + v_x)^2\}\nonumber
\end{gathered}
\end{equation}

when $p = 1$.

\

\textbf{ii.b)} {\it Furthermore, for}
\begin{equation}
\varphi = u, \quad \psi = v,\nonumber
\end{equation}

$Y_1$ also provides
\begin{equation}
\begin{gathered}
C^t = 2(u^2 + v^2), \\ C^x = 3a(u^2 + v^2)^2 + 2c[(u^2 + v^2)_{xx} - 3(u_x^2 + v_x^2)]\nonumber
\end{gathered}
\end{equation}

when $p = 2$.

\end{document}